\documentclass[12pt]{article}
\usepackage[a4paper]{geometry}
\usepackage{amssymb,amsmath,amsthm}
\usepackage{graphicx,cite,xcolor}
\usepackage[colorlinks=true,citecolor=black,linkcolor=black,urlcolor=blue]{hyperref}

\theoremstyle{plain}
\newtheorem{theorem}{Theorem}

\newtheorem{corollary}[theorem]{Corollary}

\theoremstyle{definition}

\theoremstyle{remark}

\title{A note on extremal Sombor indices\\ of trees with a given degree sequence}

\author{Ivan Damnjanovi\'c\thanks{ID is supported by Diffine LLC.}\\
\small University of Ni\v s, Faculty of Electronic Engineering, Ni\v{s}, Serbia\\[-0.4ex]
\small\tt ivan.damnjanovic@elfak.ni.ac.rs\\
\small Diffine LLC, San Diego, California, USA\\[-0.4ex]
\small\tt ivan@diffine.com
\and Marko Milo\v{s}evi\'c\\
\small University of Ni\v s, Faculty of Sciences and Mathematics, Ni\v s, Serbia\\[-0.4ex]
\small\tt markom@pmf.ni.ac.rs
\and Dragan Stevanovi\'c\thanks{DS is supported by the Serbian Ministry of Education, Science and Technological Development through the Mathematical Institute of SASA, and by the project F-159 of the Serbian Academy of Sciences and Arts.}\\
\small Mathematical Institute of the Serbian Academy of Sciences and Arts,\\[-0.4ex]
\small Belgrade, Serbia\\[-0.4ex]
\small\tt dragan\_stevanovic@mi.sanu.ac.rs
}

\begin{document}

\maketitle

\begin{abstract}
We note here that the problem of determining extremal values of Sombor index for trees with a given degree sequence fits within the framework of results by Hua Wang from [Cent. Eur. J. Math. 12 (2014) 1656--1663], implying that the greedy tree has the minimum Sombor index, while an alternating greedy tree has the maximum Sombor index.
\end{abstract}

\newpage
In a recent private communication, Ivan Gutman asked a number of colleagues 
to characterize tree(s) of order~$n$ with the given degree sequence~$\mathcal{D}$
whose Sombor indices are minimum and maximum.
Recall that Sombor index $SO(G)$ of a graph $G=(V,E)$ is defined in~\cite{gutm} as
$$
SO(G) = \sum_{uv\in E} \sqrt{d_u^2+d_v^2},
$$
where $d_u, d_v$ are degrees of the vertices $u, v\in V$.
Observing that switching edges (i.e., delete edges $ab$ and~$cd$, and add edges $ac$ and~$bd$)
decreases Sombor index under suitable conditions, while keeping degrees intact,
we quickly jumped in to show in~\cite{dast} that 
a greedy tree necessarily attains the minimum Sombor index among trees with degree sequence~$\mathcal{D}$.
Actually, the greedy tree is the unique tree that minimizes pseudo-Sombor index (see~\cite{dast} for details),
but in principle there may exist other trees with the same minimum value of Sombor index 
as the greedy tree for given~$\mathcal{D}$.
A more detailed reading of references by one of us
during the subsequent attempt to determine trees with the maximum Sombor index,
revealed that this problem actually fits within the framework of results by Hua Wang~\cite{wang},
which quickly implies both that the greedy tree has the minimum and
that an alternating greedy tree has the maximum value of Sombor index
among trees with degree sequence~$\mathcal{D}$.

For the sake of completeness, 
let us present here the result of Wang~\cite{wang}.
Assume that the degree sequence~$\mathcal{D}$ is ordered in a non-increasing order,
and denote by \linebreak $d_1\geq\dots\geq d_m$ the degrees of internal vertices
(i.e., the elements of~$\mathcal{D}$ that are greater than one).
Both the greedy tree and the alternating greedy tree are constructed algorithmically.
The greedy tree is constructed as follows:
\begin{enumerate}
\item[(g1)] Label the root with the largest degree~$d_1$;
\item[(g2)] Label the neighbors of the root as $d_2$, $d_3$, \dots,
            assigning to each next neighbor the largest available degree;
\item[(g3)] For each labelled vertex in the current level, 
            considered in a non-increasing order of labels,
            label its children in turn with the largest available degree;
\item[(g4)] Repeat (g3) as long as there are available internal degrees,
            then add necessary number of leaves so that the degree of each vertex is equal to its label.
\end{enumerate}
The alternating greedy tree is constructed by a recursive procedure:
\begin{enumerate}
\item[(a1)] If $m-1\leq d_m$, 
            the alternating greedy tree is a tree rooted at the vertex~$r$ with $d_m$ children,
            among which $d_m-m+1$ are leaves, while the remaining children have degrees $d_1,\dots,d_{m-1}$;
\item[(a2)] If $m-1\geq d_m+1$,
            create a subtree~$T$ rooted at~$r$ with $d_m-1$ children with degrees $d_1,\dots,d_{d_m-1}$;
\item[(a3)] Let $S$ be the alternating greedy tree corresponding to the sequence $(d_{d_m},\dots,\linebreak d_{m-1})$ of internal degrees,
            and let $v$ be a leaf of~$S$ with the smallest degree of its neighbor.
            The alternating greedy tree for the sequence $(d_1,\dots,d_m)$ is obtained 
            by identifying the root~$r$ of~$T$ with~$v$ in~$S$.
\end{enumerate}

\begin{figure}[ht]
\centering
\includegraphics[width=\textwidth]{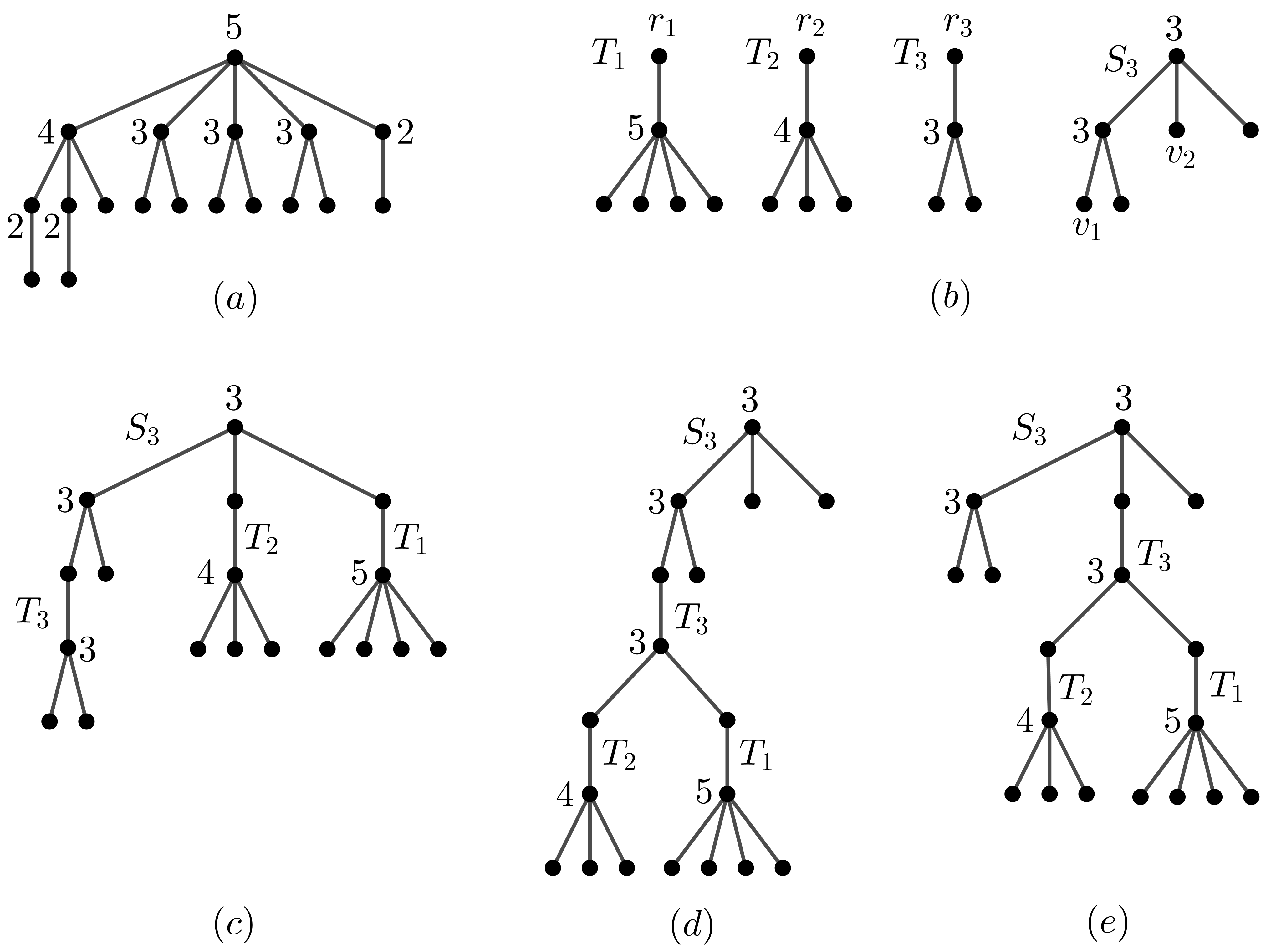}
\caption{Illustration of constructions of: (a) the greedy tree,
         (b) the constituents of alternating greedy trees, and 
         (c--e) some feasible alternating greedy trees
         for the sequence $(5,4,3,3,3,2,2,2)$ of internal vertex degrees.}
\label{fig-1}
\end{figure}

These two constructions are illustrated in Fig.~\ref{fig-1}
for the internal degree sequence $(5,4,3,3,3,2,2,2)$.
The greedy tree, shown in Fig.~\ref{fig-1}(a), is produced uniquely by the steps (g1)--(g4).
However, several non-isomorphic alternating greedy trees may be produced by the steps (a1)--(a3),
as it can happen that the leaf~$v$ in the step~(a3) can be selected in different non-isomorphic ways.
Namely, step~(a2) applied to $(5,4,3,3,3,2,2,2)$ produces the subtree~$T_1$ with the root~$r_1$,
leaving the subsequence $(4,3,3,3,2,2)$ for which one has to produce an alternating greedy tree in step~(a3).
This calls step~(a2) recursively to produce the subtree~$T_2$ with the root~$r_2$,
and leaves the subsequence $(3,3,3,2)$. 
Another recursive call to step~(a2) produces the subtree~$T_3$ with the root~$r_3$,
which leaves the subsequence~$(3,3)$ for which step~(a1) produces the alternating greedy tree~$S_3$.
All these ``constituents'' are shown in Fig.~\ref{fig-1}(b).
However, going back from these recursive calls and continuing with step~(a3)
yields several possible choices for the choice of the leaf~$v$.
First, the root~$r_3$ of~$T_3$ may be identified with either of the leaves $v_1$ and~$v_2$ of~$S_3$.
If $r_3$ is identified with~$v_1$, as done in Fig.~\ref{fig-1}(c) and \ref{fig-1}(d),
then the root~$r_2$ of~$T_2$ may be further identified with 
either one of the remaining leaves of~$S_3$ or one of the leaves of~$T_3$ in the newly formed tree.
After this is done, there are still several choices left
for the choice of the leaf which should be identified with the root~$r_1$ of~$T_1$.
Fig.~\ref{fig-1}(c)--(e) shows some of the final resulting alternating greedy trees
(where in Fig.~\ref{fig-1}(e) the root~$r_3$ was initially identified with the leaf~$v_2$).

Generalizing earlier results,
Wang~\cite{wang} considered for a tree~$T=(V,E)$
the general form of a topological index defined as 
$$
R_f(T) = \sum_{uv\in E} f(d_u, d_v),
$$
where $f\colon\mathbb{N}\times\mathbb{N}\to\mathbb{R}$ is a symmetric function.
He proved the following theorem.

\begin{theorem}[\!\!\cite{wang}]
\label{th-wang}
If the symmetric function $f\colon\mathbb{N}\times\mathbb{N}\to\mathbb{R}$ satisfies
\begin{equation}
\label{eq-condition}
f(x,a) + f(y,b) \geq f(y,a) + f(x,b)
\quad
\mbox{for all $x\geq y$ and $a\geq b$}
\end{equation}
(with strict inequality implied if both $x>y$ and $a>b$),
then $R_f(T)$ is maximized by the greedy tree and minimized by an alternating greedy tree
among trees with given degree sequence.
\end{theorem}

Let us now specifically define $f$ as 
$$
f(x,a)=-\sqrt{x^2 + a^2}
$$
so that $R_f(T)$ is actually minus Sombor index of~$T$.
In this case the condition~(\ref{eq-condition}) reads as
$$
-\sqrt{x^2+a^2} -\sqrt{y^2+b^2} \geq -\sqrt{y^2+a^2} -\sqrt{x^2+b^2}
$$
whenever $x\geq y$ and $a\geq b$,
which is equivalent to
$$
\sqrt{x^2+a^2} + \sqrt{y^2+b^2} \leq \sqrt{y^2+a^2} + \sqrt{x^2+b^2}.
$$
After squaring and rearranging, this is further equivalent to
$$
(x^2+a^2)(y^2+b^2) \leq (y^2+a^2)(x^2+b^2),
$$
and further to
$$
0 \leq (a^2-b^2)(x^2-y^2),
$$
which is certainly satisfied whenever $x\geq y$ and $a\geq b$
(with strict inequality if both $x>y$ and $a>b$).
Hence Theorem~\ref{th-wang} holds for minus Sombor index,
leading to the following corollary for Sombor index itself.

\begin{corollary}
Sombor index is minimized by the greedy tree and maximized by an alternating greedy tree
among trees with given degree sequence.
\end{corollary}


\begin{thebibliography}{0}
\bibitem{dast} I. Damnjanovi\'c, D. Stevanovi\'c,
    Greedy trees have minimum Sombor indices,
    \href{https://arxiv.org/abs/2211.05559}{\texttt{arXiv:2211.05559}} (2022).

\bibitem{gutm} I. Gutman, 
    Geometric approach to degree-based topological indices: Sombor indices,
    \textit{MATCH Commun. Math. Comput. Chem.\/} \textbf{86} (2021) 11--16.

\bibitem{wang} H. Wang,
	Functions on adjacent vertex degrees of trees with given degree sequence,
	\textit{Cent. Eur. J. Math.\/} \textbf{12} (2014) 1656--1663.

\end{thebibliography}
\end{document}